\documentclass[11pt]{article} 
\usepackage{amssymb}
\usepackage{amsmath}
\usepackage{amsthm} 
\theoremstyle{plain}
\newtheorem{thm}{\bfseries Theorem}
\newtheorem{lem}{\bfseries Lemma}
\newtheorem{prop}{\bfseries Proposition}
\newtheorem{defn}{\bfseries Definition}
\theoremstyle{definition} 
\theoremstyle{remark}
\newtheorem{rem}{\bfseries Remark}
\def\SR{{\mathsf R}}
\def\SF{{\mathsf F}}
\def\sx{{\mathsf x}}
\def\sy{{\mathsf y}}
\def\SX{{\mathsf X}}
\def\SY{{\mathsf Y}}
\def\SZ{{\mathsf Z}}
\def\ST{{\mathsf T}}
\def\SW{{\mathsf W}}
\def\SP{{\mathsf P}}

\def\SM{{\mathsf M}}
\def\SH{{\mathsf H}}
\newcommand{\dil}[1]{\left\langle #1 \right\rangle_{\!q}}
\newcommand{\ddil}[1]{{\mathbb R}_{#1}} 

\newcommand{\Dil}[1]{\biggl\langle #1\biggr\rangle_{\!q}}
\begin{document}

\begin{center}
{\bf\Large\sc 
Tetrahedron equation\\
 and cyclic quantum dilogarithm identities  
} 

\vspace*{3mm}
{\sc Andrei Bytsko \ and \ Alexander Volkov }
\end{center}

\begin{abstract}
\noindent
We establish a hierarchy of quantum dilogarithm identities 
associated to a sequence of triangular shaped quivers. 
The tetrahedron equation plays a key role in our 
construction.

\end{abstract}

\section{Introduction. Main results.}
 
Fix $q \in (0,1)$.  
The {\em quantum exponential function}
is the following formal series
\begin{equation}\label{Sx}
  \dil{x} =  \sum_{n = 0}^{\infty}  
   \frac{(-x)^n }{(1-q)\ldots(1-q^{n})} \,.
\end{equation} 
It is well known that  if $\SX$ and $\SY$ are two 
$q$--commuting indeterminates,  
that is they satisfy the commutation relation
\begin{equation}\label{qcomm}
\SY \SX = q \, \SX \SY \,,
\end{equation}  
then the following identities hold:
\begin{eqnarray}
\label{qexp} 
{}&  \dil{\SX} \,  \dil{\SY}  =  
 \dil{ \SX + \SY} \,, & \\[1mm]
{}&  \label{pent} 
 \dil{\SX} \, \dil{\SX\SY} \, \dil{\SY}  =  
 \dil{\SY} \, \dil{\SX} \,. &
\end{eqnarray}
The first one is due to Sch\"utzenberger \cite{Sch1}.
The second identity was found in \cite{FK1} and is now 
commonly called the {\em pentagon identity}.
It is also often called a {\em quantum dilogarithm identity}
because it is closely related to the five--term 
dilogarithm identity \cite{FK1,Vo1}. 
\vspace{1mm} 
 
Now take three pairwise $q$--commuting indeterminates 
$\SX$, $\SY$, and $\SZ$,
\begin{equation}\label{XYZ}
\SY \SX = q \, \SX \SY \,, \qquad
\SX \SZ = q \, \SZ \SX \,, \qquad
\SZ \SY = q \, \SY \SZ \,.
\end{equation} 
Following \cite{KV1,Se1}, we 
utilize the pentagon relation twice and find that
\begin{equation}\label{twice}  
 \dil{\SX} \,  \dil{\SZ} \, \dil{\SX\SY} \, \dil{\SY}  =  
\dil{\SZ} \, \dil{\SZ\SX} \, \dil{\SX} \, 
\dil{\SX\SY} \, \dil{\SY} = 
\dil{\SZ} \, \dil{\SZ\SX} \, \dil{\SY} \, \dil{\SX} . 
\end{equation}
Note that $\SX\SY$ and $\SZ$ commute. 
Permuting $\dil{\SX\SY}$ and $\dil{\SZ}$ on the l.h.s, 
we bring identity (\ref{twice}) to the form
$\ST = \rho (\ST)$, 
where $\ST=\dil{\SX} \, \dil{\SX\SY} \, \dil{\SZ} \, \dil{\SY}$
and $\rho: \SX \to \SZ, \SY \to \SX,  \SZ \to \SY$ 
is a third order automorphism of the associative algebra
defined by presentation~(\ref{XYZ}). 
Hence follows immediately the following triple identity 
\begin{equation}\label{444}
 \ST = \rho (\ST) = \rho \left( \rho (\ST) \right) \,,
\end{equation} 
or, explicitly,
\begin{equation}\label{44}
\dil{\SX} \, \dil{\SX\SY} \, \dil{\SZ} \, \dil{\SY}  =  
\dil{\SZ} \, \dil{\SZ\SX} \, \dil{\SY} \, \dil{\SX}  =
\dil{\SY} \, \dil{\SY\SZ} \, \dil{\SX} \, \dil{\SZ} . 
\end{equation}
In view of (\ref{444}), we will say that (\ref{44}) 
is a {\em cyclic} quantum dilogarithm identity.  
\vspace{2mm} 
 
The goal of this paper is to obtain a hierarchy of
cyclic quantum dilogarithm identities in which (\ref{44})
would be the first nontrivial member.
For this purpose we will introduce an 
algebra ${\mathcal T}_N$ with generators
assigned to the vertices of a certain quiver~$Q_N$.
In what follows, $N$ stands for an integer
number greater than one.

\begin{defn} 
 The quiver $Q_N$ is an oriented graph 
 with vertices which are labelled by pairs
of integer numbers
% are the nodes of a square lattice with coordinates 
$(i,j)$ such that $1\leq i <j \leq N$. 
The directed edges go from $(i,j)$ to $(i,j+1)$, 
{}from $(i,j)$ to $(i+1,j)$, and 
{}from $(i+1,j+1)$ to $(i,j)$. 
\end{defn}

Thus, $Q_N$ has $N \choose {2}$ vertices.
For instance, the quivers $Q_3$ and $Q_4$ are

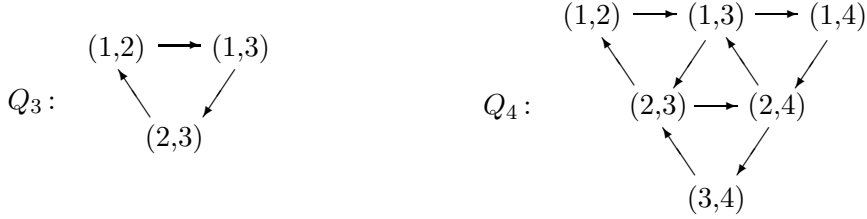
\begin{picture}(300,90)
\put(0,38){$Q_3$\,:}
\put(30,58){(1,2)}
\put(57,62){\vector(1,0){16}}
\put(77,58){(1,3)}
\put(54,35){\vector(-2,3){12}}
\put(52,24){(2,3)}
\put(86,53){\vector(-2,-3){12}}
\put(180,36){$Q_4$\,:}
\put(210,70){(1,2)}
\put(237,74){\vector(1,0){16}}
\put(257,70){(1,3)}
\put(237,47){\vector(-2,3){12}}
\put(260,13){\vector(-2,3){12}}
\put(284,47){\vector(-2,3){12}}
\put(280,36){(2,4)}
\put(283,74){\vector(1,0){16}}
\put(303,70){(1,4)}
\put(235,36){(2,3)}
\put(260,39){\vector(1,0){16}}
\put(257,1){(3,4)}
\put(288,31){\vector(-2,-3){12}}
\put(264,64){\vector(-2,-3){12}}
\put(310,64){\vector(-2,-3){12}}
\end{picture}
\vspace{1mm}

Given $Q_N$, one defines in the standard way
its skew--symmetric incidence matrix $B$:
$B_{(i,j),(i',j')}=-B_{(i',j'),(i,j)}=1$ 
if there is a directed edge  
going from $(i,j)$ to $(i',j')$, and
$B_{(i,j),(i',j')}=0$ if the
vertices $(i,j)$ and $(i',j')$ are not connected.

\begin{defn}\label{QA} 
The algebra ${\mathcal T}_N$ associated with 
the quiver $Q_N$ is a unital associative
algebra over ${\mathbb R}$ with ${N} \choose {2}$
generators $\SZ_{ij}$, $1\leq i <j \leq N$,
and the following defining relations:
\begin{align}\label{ZZ}
\SZ_{ij} \, \SZ_{i'j'}= 
 q^{ B_{(i,j),(i',j')}} \, \SZ_{i'j'} \, \SZ_{ij} \,.
\end{align} 
\end{defn}

In particular, ${\mathcal T}_2$ is generated 
by a single generator. The algebra ${\mathcal T}_3$
(with generators
$\SZ_{12}$, $\SZ_{13}$, $\SZ_{23}$ renamed
$\SX$, $\SZ$, $\SY$)
coincides with the algebra
defined by presentation~(\ref{XYZ}).

\begin{rem}\label{AAA} 
Let $I$ be a subset of vertices of $Q_N$.
Consider a linear homomorphism $\sigma_I$
such that 
$\sigma_I(\SZ_{ij})=0$ if $(i,j) \in I$
and $\sigma_I(\SZ_{ij})=\SZ_{ij}$ otherwise. 
Since relations (\ref{ZZ}) are homogeneous,
$\sigma_I$ is an algebra homomorphism from
${\mathcal T}_N$ to its subalgebra. 
In particular, sending $\SZ_{iN}$ to zero 
for all~$i$, we reduce ${\mathcal T}_N$
to ${\mathcal T}_{N-1}$. Thus, 
we  have a chain of subalgebra 
inclusions:
${\mathcal T}_2 \subset {\mathcal T}_3 \ldots
 \subset {\mathcal T}_{N-1} \subset {\mathcal T}_N$.
\end{rem}

\begin{rem}\label{remcent} 
The center ${\mathcal Z}_N$ of the algebra 
${\mathcal T}_N$ has dimension $N/2$ if
$N$ is even and $(N-1)/2$ if $N$ is odd.
(See Section 2.2 for more details.)
\end{rem}
 
The quiver $Q_N$ is mapped to itself
by a clockwise rotation about its geometric 
center  by~$2\pi/3$.
Hence, the algebra ${\mathcal T}_N$
admits a third order automorphism $\rho$,  
\begin{equation}\label{rhodef}
  \rho \bigl(\SZ_{ij} \bigr)
  = \SZ_{j-i,N+1-i} \,. 
\end{equation} 

Furthermore, the quiver $Q_N$ is mapped into
a dual quiver (were all the arrows are reverted) 
by a reflection with respect to either of
its three axes of symmetry passing through the corner vertices. 
Hence, the algebra ${\mathcal T}_N$
admits three involutive anti--automorphisms:
\begin{equation}\label{mudef}
  \mu_1 \bigl(\SZ_{ij} \bigr)
  = \SZ_{j-i,j} \, \qquad 
  \mu_2 \bigl(\SZ_{ij} \bigr)
  = \SZ_{N+1-j,N+1-i} \, \qquad 
  \mu_3 \bigl(\SZ_{ij} \bigr)
  = \SZ_{i,N+1+i-j} \,. 
\end{equation}
Note that
\begin{equation}\label{murho}
  \mu_1 \circ \mu_2 = \mu_2 \circ \mu_3 
  = \mu_3 \circ \mu_1 = \rho \,. 
\end{equation}
\vspace*{1mm}

Let us adopt the following notations for 
products of non--commuting factors.
Let $\prec$ be the lexicographic order
relation of elements of a set $\Lambda \subset {\mathbb Z}^m$.  
Then $\prod\limits_{\lambda \in \Lambda}^{\rightarrow}  f_\lambda$
and $\prod\limits_{\lambda \in \Lambda}^{\leftarrow}  f_\lambda$
stand for ordered products, where 
$f_\lambda$ is put to the right (respectively,
to the left) of all
$f_{\lambda'}$ such that $\lambda' \prec \lambda$.
In particular, if $\Lambda = [1,k] \subset {\mathbb Z}$, then
$\prod\limits_{ \lambda \in \Lambda}^{\rightarrow}  f_\lambda
= f_1 \ldots f_k$\ 
and $\prod\limits_{\lambda \in \Lambda}^{\leftarrow}  f_\lambda
= f_k \ldots f_1$.

Let $\Lambda_N \subset {\mathbb Z}^3$  be the
following discrete tetrahedron containing 
${{N+1} \choose {3}}$ points:
\begin{equation}
 \Lambda_N = 
 \bigl\{ \lambda=(a,b,c) \ | \ 
 	1 \leq a < b < c \leq N+1 \bigr\} \,.
\end{equation}
To each point $\lambda=(a,b,c) \in \Lambda_N$ we
associate the following element of
${\mathcal T}_{N}$:
\begin{equation}\label{Zkdef}  
   \ddil{ \lambda }  =
   \Dil{\,
  \prod_{0 \leq k \leq c-b-1}^{\rightarrow} \SZ_{a+k,b+k} } \,.
\end{equation}
 
Now we can define   
an analogue of the element~$\ST$ used in~(\ref{444}).
\begin{defn}\label{defT}
$\ST_N \in {\mathcal T}_{N}$ is 
the following lexicographically
ordered product: 
\begin{equation}\label{TTlex}
 \ST_N =   
  \prod_{\lambda \in \Lambda_N}^{\rightarrow}   
  \ddil{\lambda}  \,.
\end{equation}  
\end{defn}
In particular, we have 
$\ST_2=\ddil{ 123 }  = \dil{ \SZ_{12} }$ and  
$\ST_3=
\ddil{ 123 }\ddil{ 124 }\ddil{ 134 }\ddil{ 234 } = \\
 \dil{ \SZ_{12} } \dil{ \SZ_{12} \SZ_{23} }
 \dil{ \SZ_{13} }\dil{ \SZ_{23} }$.
\vspace*{1mm}

\begin{rem} 
Some factors in the product (\ref{TTlex}) mutually commute,
as, e.g., $\ddil{ 124 }$ and $\ddil{ 134 }$ in $\ST_3$.
Therefore, $\ST_N$ admits a number of equivalent expressions 
obtained by permutations of such factors.
Some of these expressions are given in Lemma~\ref{TW23}.  
\end{rem}
\vspace{1mm}

We will say that a quantum dilogarithm identity 
is of the type $n \sim m$ if it involves $n$ quantum 
exponentials with monomial arguments on one side and $m$ on
the other. For instance, (\ref{pent}) and (\ref{44})
are of the type $2 \sim 3$ and $4 \sim 4$, respectively. 
The main result of the present work
is a family of cyclic quantum dilogarithm identities
of the type ${{N+1} \choose {3}} \sim {{N+1} \choose {3}}$.
 
\begin{thm}\label{main}
Identities
\begin{eqnarray}
\label{Tmu}
{}& 
 \ST_N = \mu_1  \bigl( \ST_N \bigr) =
  \mu_2  \bigl( \ST_N \bigr) =
  \mu_3  \bigl( \ST_N \bigr) \,, & \\[2mm]
\label{Tr}
{}&
  \ST_N =  \rho \bigl( \ST_N \bigr)
  = \rho \bigl( \rho \bigl( \ST_N \bigr) \bigr) \,, &
\end{eqnarray}
hold for any integer $N \geq 2$. 
\end{thm}

\begin{rem}
Identities (\ref{Tmu}) involve the
anti--automorphism transformations $\mu_a$
corresponding to a reflection rather than
a rotation symmetry of the quiver.
Nevertheless, two of these identities
are almost cyclic in the following sense.
We will see below (cf. Lemma~\ref{mu2T}) 
that $\mu_2$ applied to $\ST_N$ acts 
almost as an identical transformation just  
permuting some commuting quantum exponentials. 
(As an example, apply to (\ref{44}) 
an anti--automorphism which maps $\SZ$ to itself and
exchanges $\SX$ and $\SY$.)
Along with (\ref{murho}) this implies 
that, again up to a permutation of commuting factors,
$\mu_1$ applied to $\ST_N$ acts as $\rho$,
and $\mu_3$ applied to $\ST_N$ acts as 
$\rho^{-1}=\rho \circ \rho$.
\end{rem}

It is clear from the definition (\ref{Sx}) 
of the quantum exponential 
 that $\dil{0}=1$. This, along with 
Remark~\ref{AAA}, implies the following.

\begin{prop}\label{QQ}
Let $I$ be a subset of the set of vertices of $Q_N$. 
In the identity (\ref{Tr}),
replace with unity every quantum exponential 
that contains at least one $\SZ_{ij}$ with
$(i,j) \in I$. The result is 
a correct quantum dilogarithm identity.
\end{prop}

Note that the resulting reduced identity is not
necessarily cyclic. For instance, sending
$\SZ$ to $0$ in the $4 \sim 4$ identity (\ref{44}), we
obtain the pentagon identity~(\ref{pent}).
\vspace{2mm}

Quantum dilogarithm identities associated with
various quivers have close connections with 
(quantum) cluster algebras and Y--systems, 
see, e.g.~\cite{BZ1,FG1,KN1,Ke1}. 
Our motivation for studying a specific 
family of such identities is an observation that
they are {\em accessible} from 
the $4 \sim 4$ identities~(\ref{44}).
For instance, let us demonstrate how to derive the equality 
$\ST_4 =  \rho \bigl( \ST_4 \bigr)$ 
without invoking 
the pentagon relation but using only the $4 \sim 4$ identities~(\ref{44}): 
\begin{align*}
\ST_4
{}&= \dil{ \SZ_{12} } \dil{ \SZ_{12} \SZ_{23} } 
\dil{ \SZ_{12} \SZ_{23} \SZ_{34} }
\dil{ \SZ_{13} }  \dil{ \SZ_{13} \SZ_{24} }   
 \dil{ \SZ_{14} } \dil{ \SZ_{23} } 
\dil{ \SZ_{23} \SZ_{34} }  
\dil{ \SZ_{24} } \dil{ \SZ_{34} }  \\ 
{}&=  \underline{ \dil{ \SZ_{12} } \dil{ \SZ_{12} \SZ_{23} } 
\dil{ \SZ_{13} }  \dil{ \SZ_{23} }  }
\dil{ \SZ_{12} \SZ_{23} \SZ_{34} } 
\dil{ \SZ_{13} \SZ_{24} } 
\dil{ \SZ_{23} \SZ_{34} }  \dil{ \SZ_{14} }
\dil{ \SZ_{24} } \dil{ \SZ_{34} }  \\
{}& = \dil{ \SZ_{13} }  
 \dil{ \SZ_{13}\SZ_{12} } \dil{ \SZ_{23} } 
\underline{   \dil{ \SZ_{12} } 
\dil{ \SZ_{12} \SZ_{23} \SZ_{34} } 
\dil{ \SZ_{13} \SZ_{24} }
\dil{ \SZ_{23} \SZ_{34} } }
  \dil{ \SZ_{14} }
\dil{ \SZ_{24} } \dil{ \SZ_{34} } \\
{}&= \dil{ \SZ_{13} }  \dil{ \SZ_{23} } 
 \dil{ \SZ_{13}\SZ_{12} }  
\dil{ \SZ_{13} \SZ_{24} } \dil{ \SZ_{13} \SZ_{24} \SZ_{12} }
\dil{ \SZ_{23} \SZ_{34} }  \dil{ \SZ_{12} } 
 \dil{ \SZ_{14} }
\dil{ \SZ_{24} } \dil{ \SZ_{34} } \\ 
{}& =\dil{ \SZ_{13} }  \dil{ \SZ_{23} }
\dil{ \SZ_{13} \SZ_{24} } 
\underline{ \dil{ \SZ_{13}\SZ_{12} } 
\dil{ \SZ_{13}  \SZ_{12} \SZ_{24} }
\dil{ \SZ_{14} } \dil{ \SZ_{24} } }
\dil{ \SZ_{23} \SZ_{34} }  \dil{ \SZ_{34} } \dil{ \SZ_{12} } \\
{}&= \dil{ \SZ_{13} }  \dil{ \SZ_{23} }
\dil{ \SZ_{13} \SZ_{24} } 
 \dil{ \SZ_{14} } 
\dil{ \SZ_{14} \SZ_{13} \SZ_{12}} 
\dil{ \SZ_{24} } \dil{ \SZ_{13}\SZ_{12} }  
\dil{ \SZ_{23} \SZ_{34} }  \dil{ \SZ_{34} } \dil{ \SZ_{12} } \\
{}&= \dil{ \SZ_{13} } \dil{ \SZ_{13} \SZ_{24} }  
\dil{ \SZ_{14} } 
\dil{ \SZ_{14} \SZ_{13} \SZ_{12}} 
\underline{ \dil{ \SZ_{23} }
  \dil{ \SZ_{23} \SZ_{34} }  \dil{ \SZ_{24} }
\dil{ \SZ_{34} } }
 \dil{ \SZ_{13}\SZ_{12} } \dil{ \SZ_{12} }
\\
{}&= 
\underline{ \dil{ \SZ_{13} } \dil{ \SZ_{13} \SZ_{24} }  
\dil{ \SZ_{14} } \dil{ \SZ_{24} } }
\dil{ \SZ_{14} \SZ_{13} \SZ_{12}}     
 \dil{ \SZ_{24} \SZ_{23} }
\dil{ \SZ_{34} } \dil{ \SZ_{23} }
 \dil{ \SZ_{13}\SZ_{12} }   \dil{ \SZ_{12} }
\\ 
{}&= 
\dil{ \SZ_{14} } 
\dil{ \SZ_{14} \SZ_{13} } \dil{ \SZ_{24} }
 \dil{ \SZ_{13} }  
\dil{ \SZ_{14} \SZ_{13} \SZ_{12}}     
 \dil{ \SZ_{24} \SZ_{23} }
\dil{ \SZ_{34} }  
 \dil{ \SZ_{13}\SZ_{12} } \dil{ \SZ_{23} } \dil{ \SZ_{12} }
\\ 
{}& = \dil{ \SZ_{14} }
\dil{ \SZ_{14} \SZ_{13} }  
\dil{ \SZ_{14} \SZ_{13} \SZ_{12} }  
\dil{ \SZ_{24} } \dil{ \SZ_{24} \SZ_{23} }
\dil{ \SZ_{34} }  \dil{ \SZ_{13} }   
 \dil{ \SZ_{13}\SZ_{12} } \dil{ \SZ_{23} }
\dil{ \SZ_{12} } 
\\{}& = \rho \bigl(\ST_4 \bigr).
\end{align*}  
The underlined terms were transformed by
applying the $4\sim 4$ identities~(\ref{44}).
The remaining transformations changed only the
order of commuting factors or the order of
commuting generators in the arguments of 
quantum exponentials. 
\vspace{1mm}

We will give a proof of the identities listed in
Theorem~\ref{main} which makes it evident 
that they are accessible from the $4 \sim 4$ 
identities~(\ref{44}) for all~$N$.
The origin of this accessibility is that 
these identities  
stem from identities for certain words in
a group whose generators satisfy 
the {\em tetrahedron equation},
\begin{equation}\label{RTE0}  
 R_{abc} \, R_{abd} \, R_{acd} \, R_{bcd} =
  R_{bcd} \, R_{acd} \, R_{abd} \, R_{abc} \,.
\end{equation} 
\vspace{1mm}

The paper is organized as follows. Section~2
contains auxiliary statements which we 
need to combine together in order to prove 
Theorem~\ref{main}. Namely, in Section~2.1
we consider families of groups ${\mathcal B}(n,N)$ whose 
generators satisfy the Yang--Baxter equation
(for $n=2$), the tetrahedron equation (for $n=3$),
or their higher analogues (for \hbox{$n\geq 4$}). 
The two key technical
results here are an identity for certain words
containing all the generators of the group ${\mathcal B}(n,N)$
and relation of these words to the 
element $\ST_N$.
In Section~2.2 we describe the center of ${\mathcal T}_N$.
In Section~2.3 we consider a local tensor space 
representation $\phi$ for ${\mathcal B}(3,N)$.
In Section~2.4 we explain how evaluation of the above 
mentioned identity for words of 
${\mathcal B}(3,N)$ in the representation $\phi$ yields
ultimately the desired quantum dilogarithm
identities for an arbitrary~$N$. Appendix contains
proofs of all statements given in Section~2.

\section{Main technical ingredients}

Below we assume that $N$ and $n$ are positive integers
 and $N \geq n$.
 
\subsection{A group with $n$--simplex relations}

\begin{defn}\label{BRN}
${\mathcal B}(n,N)$ is a group with
${N} \choose {n}$ generators $R_{a_1,\ldots,a_n}$,
where 
$1 \leq a_1 < a_2 < \ldots < a_n \leq N$.
The group is defined by the following presentation:
\begin{itemize}
\item
The generators commute,
\begin{equation}\label{Rcomm}
  R_{a_1,\ldots,a_n} \, R_{b_1,\ldots,b_n} =
  R_{b_1,\ldots,b_n} \, R_{a_1,\ldots,a_n}  \,,
\end{equation}
unless the set $\{a_1,\ldots,a_n\} \bigcap \{b_1,\ldots,b_n\}$ 
contains exactly $(n-1)$ element. 
\item
If $N > n$, 
the generators 
satisfy the following ${N} \choose {n+1}$ relations:
\begin{equation}\label{Rgen}
     \prod_{1\leq j \leq n+1}^{\rightarrow} 
   R_{a_1,\ldots,\check{a}_j,\ldots,a_{n+1}}
   = \prod_{1\leq j \leq n+1}^{\leftarrow} 
   R_{a_1,\ldots,\check{a}_j,\ldots,a_{n+1}} \,,
\end{equation}
where $\check{a}_j$ is dropped. 
\end{itemize}
\end{defn}

For $n=1$, relations (\ref{Rgen}) 
imply commutativity,
$R_a R_b = R_b R_a$, so that ${\mathcal B}(1,N)$
is an abelian group with $N$ generators.

For $n=2$, relations (\ref{Rgen}) have the form
of the {\em Yang--Baxter equation},
\begin{equation}\label{RYB}  
 R_{ab} \, R_{ac} \, R_{bc} = R_{bc} \, R_{ac} \, R_{ab} \,,
 \qquad a< b < c \,,
\end{equation}
and $R_{ab}$ commutes with $R_{a'b'}$ if they have no
common index.

For $n=3$, relations (\ref{Rgen}) have the form
of the {\em tetrahedron equation},
\begin{equation}\label{RTE}  
 R_{abc} \, R_{abd} \, R_{acd} \, R_{bcd} =
  R_{bcd} \, R_{acd} \, R_{abd} \, R_{abc} \,,
 \qquad a< b < c < d \,,
\end{equation}
and $R_{abc}$ commutes with $R_{a'b'c'}$ 
unless they have exactly two
common indices.

Recall that, given a set $\Lambda \subset {\mathbb Z}^m$, 
we use the symbol 
$\prod\limits^{\rightarrow}_{\lambda\in\Lambda} f_\lambda$ 
to denote the lexicographically ordered product of 
non--commuting factors. 

\begin{defn}\label{defWW}
The word $W(n,N) \in {\mathcal B}(n,N)$ is
the lexicographically ordered product of
all generators of ${\mathcal B}(n,N)$, that is
\begin{equation}\label{lex}  
  W(n,N) = \prod_{1\leq a_1 <a_2 \ldots  < a_n \leq N 
  }^{\rightarrow}  
  R_{a_1,\ldots,a_{n}}   \,.
\end{equation} 
\end{defn}
In particular, we have $W(n,n) = R_{1,\ldots,n}$ for all $n$ 
and $W(1,N) = R_1 \ldots R_N$ for all~$N$.
\vspace*{2mm}

Let $\prec'$ be the colexicographic order
relation of elements of a set $\Lambda \subset {\mathbb Z}^m$.  
That is, components of elements of $\Lambda$ are compared 
starting from the right.  
For instance, $(a,b) \prec' (c,d)$ iff $b < d$ or $b=d$ and $a<c$.
We will denote by   
$\mathop{{\prod}'}\limits^{\rightarrow}_{\lambda \in \Lambda}  f_\lambda$
and 
$\mathop{{\prod}'}\limits^{\leftarrow}_{\lambda \in \Lambda}  f_\lambda$
ordered products where 
$f_\lambda$ is put to the right (respectively,
to the left) of all
$f_{\lambda'}$ such that $\lambda' \prec' \lambda$.

\begin{lem}\label{defWalt}
Define
\begin{equation}\label{Wexp}  
  W'(n,N) = 
 \mathop{{\prod}'}\limits^{\rightarrow}_{1\leq a_1 <a_2 \ldots  < a_n \leq N }    
  R_{a_1,\ldots,a_{n}}   \,.
\end{equation} 
Then we have the equality
\begin{equation}\label{WWp}  
  W(n,N) = W'(n,N)   \,.
\end{equation} 
\end{lem}
For instance, (\ref{lex}) yields 
$W(2,4)  = R_{12} R_{13} R_{14} R_{23} R_{24} R_{34}$
while (\ref{Wexp}) yields
 $W'(2,4)  = R_{12} R_{13} R_{23} R_{14} R_{24} R_{34}$.
These two words coincide since, by (\ref{Rcomm}),
$R_{14}$ and $R_{23}$ commute.
It is also true in general that one needs to use only
the commutativity relations (\ref{Rcomm}) in order
to change the order of factors in (\ref{Wexp})
to match that in~(\ref{lex}).
\vspace{2mm} 

Let $\eta$ be an involutive anti--automorphism
of ${\mathcal B}(n,N)$ such that
\begin{equation}\label{eta}  
   \eta \bigl( R_{a_1,\ldots,a_{n}} \bigr)
  =  R_{a_1,\ldots,a_{n}} 
\end{equation}
for all the generators of~${\mathcal B}(n,N)$.
Define 
\begin{equation}\label{bomW}  
  \bar{W}(n,N) = \eta \bigl( W(n,N) \bigr)    \,,
\end{equation}
which is a word with the order of factors reverse 
to that of~$W(n,N)$. 
Our first key technical statement is the following.
\begin{thm}\label{multitetra}
For all $N \geq n$, we have the equality
\begin{equation}\label{WWr}
 W(n,N)  = \bar{W}(n,N)   \,.
\end{equation}
\end{thm}

\begin{rem}
A simple inspection of the proof given in the Appendix shows 
that the invertibility of $R$'s is not really needed. 
Thus, Theorem~\ref{multitetra} holds also
if ${\mathcal B}(n,N)$ is a semigroup. 
\end{rem}

\begin{rem}\label{SYZ}
The presented proof provides a constructive recursive 
procedure that transforms $W(n,N)$ into $\bar{W}(n,N)$.
The transformation involves trivial moves 
based on (\ref{Rcomm}) and moves
${\mathcal R}_{a_1,\ldots,a_{n+1}}$ which 
transform the l.h.s. of (\ref{Rgen}) into
its r.h.s.  For instance,
${\mathcal R}_{123}\bigl( W(2,4)\bigr)=
R_{23} R_{13} R_{12} R_{14} R_{24} R_{34}$. 
Let ${\mathcal W}(n,N)$ and
$\bar{\mathcal W}(n,N)$ stand for the compositions of 
such moves in which $\mathcal R$'s are ordered in the 
same way as $R$'s are ordered in $W(n,N)$
and $\bar{W}(n,N)$. A simple inspection of the proof
shows that the l.h.s. of (\ref{WWr}) is 
transformed into its r.h.s. by  
${\mathcal W}^*(n+1,N)$, where star means that 
trivial moves are included 
when necessary. Moreover, the same transformation 
is achieved by $\bar{\mathcal W}^*(n+1,N)$
if we start each reordering not from the left but from
the right. Thus, the moves satisfy
the identity 
${\mathcal W}^*(n,N) = \bar{\mathcal W}^*(n,N)$
analogous to~(\ref{WWr}).
\end{rem}
\vspace*{1mm}

Eqs. (\ref{lex}) and (\ref{Wexp}) for $n=2$
and $n=3$ read 
% the words $W(2,N)$ and $W(3,N+1)$ are given by
\begin{align}
\label{W2lex}  
&  W(2,N) = \prod_{1\leq a < b \leq N}^{\rightarrow}   
  R_{ab}  , \qquad\quad\ 
  W'(2,N)= 
\mathop{{\prod}'}\limits^{\rightarrow}_{1\leq a< b \leq N }    
  R_{ab} \,,  \\
\label{W3lex}
&  W(3,N) = \prod_{1\leq a < b < c \leq N }^{\rightarrow}   
  R_{abc}  , \qquad  
  W'(3,N) = 
  \mathop{{\prod}'}\limits^{\rightarrow}_{1\leq a< b <c \leq N }  
  R_{abc} \,.
\end{align} 
Recall that $\ddil{abc}$ was defined in (\ref{Zkdef}).
We introduce also the following elements of~${\mathcal T}_N$:
\begin{equation}\label{Rkdef}   
  \ddil{  ab  }  \, = \!
  \prod_{b+1 \leq c \leq N+1}^{\rightarrow} \ddil{ abc } \,.
\end{equation}

Let us introduce the following homomorphisms from
${\mathcal B}(2,N)$ and ${\mathcal B}(3,N+1)$
to ${\mathcal T}_N$  (it should be stressed 
that they are not algebra homomorphisms):
\begin{align}
\label{RabZ}  
{} \sharp :  R_{ab} \to \ddil{ ab } \,, \quad
{}   \star :  R_{abc} \to \ddil{ abc } \,, \quad
  \star\star :  R_{abc} \to \ddil{ a,c+a-b,c } \,.
\end{align}

\begin{lem}\label{TW23} 
The element $\ST_N$ defined in~(\ref{TTlex})
can be obtained as follows:
\begin{align}
\label{TWRZa}  
{}& \ST_N = \bigl( W(2,N) \bigr)^\sharp 
 = \bigl( W'(2,N) \bigr)^\sharp \,,\\ 
\label{TWRZb}  
{}& \ST_N = \bigl( W(3,N+1) \bigr)^\star = 
\bigl( W'(3,N+1) \bigr)^{\star} \,, \\
\label{TWRZc}  
{}& \ST_N = \bigl( W(3,N+1) \bigr)^{\star\star} = 
\bigl( W'(3,N+1) \bigr)^{\star\star} \,.
 \end{align}
\end{lem}
This observation suggests that  
identities~(\ref{WWr}) can be recast in
some quantum dilogarithm identities. 
Furthermore, Remark \ref{SYZ} indicates that one needs
to use the $4 \sim 4$ relation~(\ref{44})
$N+1 \choose 4$ times in order to derive 
identity (\ref{Tr}) for given~$N$.

\subsection{The center of ${\mathcal T}_N$}
Let us assign to each vertex $(i,j)$ of the quiver $Q_N$
a non--negative integer weight $\alpha_{ij}$. 
Let $\alpha$ denote the vector
comprised of those weights taken in the lexicographic
order. For instance, 
$\alpha=(\alpha_{12},\alpha_{13},\alpha_{23})$ if $N=3$. 
  
Monomials
\begin{equation}\label{mona} 
 \SM(\alpha) =   \prod_{1\leq i < j \leq N}^{\rightarrow}   
  \SZ_{ij}^{\alpha_{ij}} \,, \qquad 
  \alpha \in {\mathbb Z}_{\geq 0}^{N(N-1)/2}
\end{equation} 
constitute a basis of ${\mathcal T}_N$.
We will say that $\alpha_{ij}$ are the weights of $ \SM(\alpha)$.
\vspace{2mm}

The algebra ${\mathcal T}_N$ has a nontrivial center ${\mathcal Z}_N$. 
In particular, for any~$N$, if $\alpha^0$ is a vector such that $\alpha^0_{ij}=1$,
then $\SM(\alpha^0)$ belongs to ${\mathcal Z}_N$.

\begin{thm}\label{center}  
Let $\chi(N)$ stand for $N/2$ if
$N$ is even and for $(N-1)/2$ if $N$ is odd.

\begin{itemize}
\item
If\ $ \SM(\alpha) \in {\mathcal Z}_N$, then
its weights 
have the following symmetries 
\begin{equation}\label{alsym}
\alpha_{ij}   = \alpha_{j-i,N+1-i} =
 \alpha_{j-i,j} =
 \alpha_{N+1-j,N+1-i} = \alpha_{i,N+1+i-j}  \,,
\end{equation}
and $\SM(\alpha)$ is 
invariant under the action of 
$\rho$ and $\mu_k$ defined in
(\ref{rhodef}) and~(\ref{mudef}).
\item
If\  $ \SM(\alpha), \SM(\alpha') \in {\mathcal Z}_N$, and
$\alpha_{1j}=\alpha'_{1j}$ for $j=2,\ldots,\chi(N) {+} 1$,
then $\alpha=\alpha'$ and hence $ \SM(\alpha)=\SM(\alpha')$.
\item
Given $N$ and an arbitrary integer sequence
$\beta_{1}, \ldots  ,\beta_{\chi(N) }$,
there exist $m \in {\mathbb Z}_{\geq 0}$ and 
a vector $\alpha \in {\mathbb Z}_{\geq 0}^{N(N-1)/2}$ such that
$\SM(\alpha) \in {\mathcal Z}_N$  and 
$\alpha_{1j}=\beta_{j-1} +m$ for $j=2,\ldots,\chi(N) {+} 1$.
\end{itemize}
\end{thm}

In other words, each monomial central element 
$ \SM(\alpha) \in {\mathcal Z}_N$
is uniquely determined by its weights
assigned to the vertices comprising a half of a boundary 
side of the quiver $Q_N$.  It follows also that 
the dimension of ${\mathcal Z}_N$ is~$\chi(N)$.

For instance, for $N=3$ we have $\chi(3)=1$ and
$ \SM(\alpha) \in {\mathcal Z}_3$ 
iff $ \alpha_{13}=\alpha_{23}=\alpha_{12}$.
Thus, ${\mathcal Z}_3$ is generated by 
$\SZ_{12} \SZ_{13} \SZ_{23} $.
For $N=4$ we have $\chi(4)=2$ and 
$ \SM(\alpha)\in {\mathcal Z}_4$ iff
$ \alpha_{14}=\alpha_{34}=\alpha_{12}$ and
$ \alpha_{23}=\alpha_{24}=\alpha_{13}$.
Thus, ${\mathcal Z}_4$ is generated by 
$\SZ_{12} \SZ_{14} \SZ_{34} $ and 
$\SZ_{13} \SZ_{23} \SZ_{24} $.

For $N=5$ we have $\chi(5)=2$ and 
$ \SM(\alpha)\in {\mathcal Z}_5$ iff
$ \alpha_{15}=\alpha_{45}=\alpha_{12}$,
$ \alpha_{14}=\alpha_{25}=\alpha_{35}=
\alpha_{23}=\alpha_{34}=\alpha_{13}$,
and $\alpha_{24}=2\alpha_{13} - \alpha_{12}$.
The last relation requires that
$2\alpha_{13} \geq \alpha_{12}$.

\begin{rem}
The last part of Theorem~\ref{center}
shows that the $\chi(N)$ weights that
define a central monomial can be 
taken almost arbitrary (possibly, up to a
 total shift by an integer $m$).
A direct inspection up to $N=8$ suggests that  
in order to have $m=0$ it suffices 
to take a non--decreasing integer sequence, 
$0 \leq \beta_{1} \leq \ldots  \leq \beta_{\chi(N) }$. 
\end{rem}

\begin{rem}
The generators $\SZ_{ij}$ of ${\mathcal T}_N$
can be constructed as exponential functions of 
canonically conjugate variables $p_k, x_k$ 
such that $[x_k, p_{k'}]=\sqrt{-1}\,\delta_{k k'}$.
Theorem~\ref{center} implies that 
the number of {\em degrees of freedom}
of ${\mathcal T}_N$, i.e., the minimal number of 
such pairs $(p_k, x_k)$ is $N(N-2)/4$ if
$N$ is even and $(N-1)^2/4$ if $N$ is odd.
\end{rem}

\subsection{A local tensor space representation for 
tetrahedron equation}

Let $V$  be a vector space. For $N \geq 2$, we 
define $S_N = V^{\otimes N(N-1)/2}$. 
The tensor components of $S_N$ will be 
labelled in the lexicographic order by pairs 
of integers $(i,j)$, where $1 \leq i < j \leq N$.
For instance, $S_3 = V_{12} \otimes V_{13} \otimes V_{23}$,
where all $V_{ij}$ are isomorphic to~$V$.

For a given $\SR \in \text{End}\,\bigl(S_3\bigr)$, 
we denote by
$\SR_{abc}$ its canonical extension to a linear 
operator on $S_N$ which acts non--trivially only 
on the tensor components $V_{ab}$, $V_{ac}$, and $V_{bc}$.
E.g., 
$\SR_{123}=\SR \otimes id \otimes id \otimes id$\;
if~$N=4$.

\begin{defn}
Let $\SR \in \text{End}\,\bigl(S_3\bigr)$ be invertible 
and let its action canonically extended 
to $S_4$ satisfy the tetrahedron equation, 
\begin{equation}\label{RRR}
  \SR_{123} \, \SR_{124} \, \SR_{134} \, \SR_{234} =
  \SR_{234} \, \SR_{134} \, \SR_{124} \, \SR_{123} \,.
\end{equation}
A {\em local tensor space representation} of 
${\mathcal B}(3,N)$ is a homomorphism 
${\mathcal B}(3,N) \to \text{End}\,\bigl(S_{N}\bigr)$ 
sending $R_{abc}$ to $\SR_{abc}$.
\end{defn}

Note that the commutativity relations (\ref{Rcomm}) hold by
construction. Indeed, if $\SR_{abc}$ and $\SR_{a'b'c'}$
have no common pair of indices then they act non--trivially
on different tensor components of~$S_N$ and hence they commute. 

\begin{rem}
An equation formally identical to (\ref{RRR}) was
considered in~\cite{FM1}. There, however, it was
treated as an operator equation on $V^{\otimes 4}$
with the identification $\SR_{123} =\SR \otimes id$ etc.
We treat (\ref{RRR}) as an equation on 
$V^{\otimes 6}$ which is the standard Zamolodchikov's
tetrahedron equation \cite{Za1} but with a non--standard
(double index) labelling of the tensor components.
\end{rem}

Let $V$ be the vector space of formal series 
in $x,x^{-1}$. Then  $S_N$ is
the vector space of formal series in  
$x_{ij}, x_{ij}^{-1}$, where $1 \leq i < j \leq N$. 

Define operators
$\sx_{ij}, \sy_{ij} \in \text{End}\,\bigl(S_N\bigr)$,
$1 \leq i < j \leq N$ such that
\begin{equation}\label{xyf}
\begin{aligned}
  \bigl(\sx_{ij} \, f)(x_{12},\ldots,x_{ij},\ldots)
 {} &=  x_{ij} \, f(x_{12},\ldots,x_{ij},\ldots) \,, \\
  \bigl(\sy_{ij} \, f)(x_{12},\ldots,x_{ij},\ldots)
 {}  &= f(x_{12},\ldots,q x_{ij},\ldots) \,,
\end{aligned}
\end{equation}
for any $f \in S_N$. These operators comprise $N \choose 2$  
$q$--commuting pairs,  
\begin{equation}\label{qxy}
 \sx_{ij} \, \sx_{i'j'} =   
  \sx_{i'j'} \, \sx_{ij} \,, \qquad
  \sy_{ij} \, \sx_{i'j'} = q^{\delta_{ii'} \delta_{jj'}} \,
  \sx_{i'j'} \, \sy_{ij} \,, \qquad
  \sy_{ij} \, \sy_{i'j'} =  
  \sy_{i'j'} \, \sy_{ij} \,.
\end{equation}

Consider  $\SF\in \text{End}\,\bigl(S_3\bigr)$
whose action on monomials is given by
\begin{equation}\label{Fdef}
  \SF : x_{12}^k \, x_{13}^l \, x_{23}^m \to
   x_{12}^k \, x_{13}^{m+k} \, x_{23}^{l-k} \,,
\end{equation}
or, equivalently, $(\SF f)(x_{12},x_{13},x_{23})
= f(\frac{x_{12}x_{13}}{x_{23}},x_{23},x_{13})$.
It is easy to check that
\begin{align}
\label{Fxx}    
{}&  \SF \, \sx_{12} = \sx_{12} \sx_{13} \sx_{23}^{-1} \SF \,, \qquad
\SF \, \sx_{13} = \sx_{23} \, \SF \,, \qquad
\SF \, \sx_{23} = \sx_{13} \, \SF \,, \\
\label{Fyy}     
{}&  \SF \, \sy_{12} = \sy_{12} \, \SF \,, \qquad
 \SF \, \sy_{13} = \sy_{12} \sy_{23} \, \SF \,, \qquad
 \SF \, \sy_{23} = \sy_{12}^{-1} \sy_{13} \, \SF \,.
\end{align}
Clearly, $\SF$ is invertible and $\SF^2 = id$.

Let
$\SR \in \text{End}\,\bigl(S_3\bigr)$
be the following operator
\begin{equation}\label{Rdef}    
  \SR = \SF \cdot
 \dil{ q^{1+\gamma} \sx_{12}   \sx_{23}^{-1} 
  \sy_{12}^{\gamma} \sy_{13}^{-\gamma} \sy_{23}^{-1-\gamma} }   \,,
\end{equation}
where $\gamma\in \mathbb Z$.  
A reader familiar with the Yang--Baxter
equation will see in (\ref{Rdef}) an analogy with the
standard ansatz for $R$--matrix,  $\SR=\SP\cdot \check{\SR}$,
where $\SP$ is the permutation, $\SP^2=id$. 

\begin{lem}\label{RSR}
Let $\SF_{abc}$ and $\SR_{abc}$ stand for the
canonical extensions of $\SF$ defined by (\ref{Fdef})
and $\SR$ defined by (\ref{Rdef}) 
to operators on $S_N$.
The homomorphisms $\theta, \phi_\gamma :
{\mathcal B}(3,N) \to \text{End}\,\bigl(S_{N}\bigr)$ 
such that $\theta(R_{abc})=\SF_{abc}$ and
$\phi_\gamma(R_{abc})=\SR_{abc}$ are 
local tensor space representations of~${\mathcal B}(3,N)$. 
\end{lem}

\begin{rem}
This Lemma generalizes two previously known 
solutions of the tetrahedron equation. Namely,
the case corresponding to $\gamma=0$
was found in~\cite{Se1}, and the case corresponding 
to $\gamma=-1$ was considered in~\cite{KV1}. 
Our proof follows closely that given in~\cite{KV1}. 
\end{rem}

\begin{lem}\label{Zrep}
Let $\sx_{ab}$,  $\sy_{ab}$,
$1 \,{\leq}\, a \,{<}\, b \,{\leq}\, N+1$
act on $S_{N+1}$ as defined in~(\ref{xyf})
and let $\sy_{ab} \equiv 1$ if $a=b$.
Then, for all $N \geq 2$,  the linear homomorphism 
$\tau : {\mathcal T}_{N} \to \text{End}\,\bigl(S_{N+1}\bigr)$
such that 
\begin{equation}\label{Zxy}
  \tau(\SZ_{ab} ) =
  q\, \frac{\sx_{ab} }{ \sx_{a+1,b+1} }
   \frac{ \sy_{a+1,b}}{\sy_{a+1,b+1}} \,
  \qquad \text{for}\quad
  1 \,{\leq}\, a \,{<}\, b \,{\leq}\, N \,,
\end{equation} 
is a faithful representation of~${\mathcal T}_{N}$.
\end{lem}

\subsection{Proof of Theorem~\ref{main}} 
One the three symmetries of $\ST_N$
presented in (\ref{Tmu}) can be established
by reordering commuting factors. 
Namely, applying  $\mu_2$ to
relations (\ref{TWRZb}) and (\ref{TWRZc})
and using Lemma~\ref{mu2T} (see Appendix~A.4), 
we obtain the following statement.
\begin{lem}\label{Tmu2} 
For all $N \geq 2$, we have
\begin{equation}\label{tmu2}
  \mu_2 \bigl( \ST_N \bigr) =   \ST_N \,.
\end{equation} 
\end{lem}
Clearly, in order to establish an analogous 
equality involving $\mu_1$ or~$\mu_3$ we 
will have to use the $4 \sim 4$ relation~(\ref{44}).
To do it for $\mu_1$, we establish a connection 
between the element $\ST_N$ evaluated in the
representation~$\tau$ and 
words from ${\mathcal B}(3,N+1)$ evaluated
in the representation $\phi_\gamma$ (we will
take $\gamma=0$ for simplicity).
\begin{lem}\label{RRR'}
For all $N \geq 2$, we have the equalities
\begin{align}
\label{RFFt}
{}&  \phi_0 \bigl(W'(3,N+1)\bigr) = 
  \theta \bigl(W'(3,N+1)\bigr) \,
  \tau \bigl( \ST_N \bigr) \,, \\[1mm]
\label{RFFtr} 
{}&   \phi_0 \bigl(\bar{W}'(3,N+1)\bigr) = 
  \theta \bigl(\bar{W}'(3,N+1)\bigr) \,
  \tau \bigl( \mu_1 (\ST_N) \bigr) \,.
\end{align} 
\end{lem}
Invoking Theorem~\ref{multitetra}, we infer that
$
\tau\bigl( \ST_N\bigr) = 
\tau \bigl( \mu_1 \bigl(\ST_N \bigr) \bigr) 
$. 
Which, by virtue of Lemma~\ref{Zrep},
implies that 
\begin{equation}\label{tmu1}
 \mu_1 \bigl(\ST_N \bigr) = \ST_N \,.
\end{equation}
Now, combining (\ref{tmu2}) with (\ref{tmu1})
and using (\ref{murho}), we conclude that
$\rho \bigl(\ST_N \bigr) = 
 (\mu_1 \circ \mu_2 )\bigl(\ST_N \bigr) = \ST_N$.
And finally, using (\ref{murho})  again, we obtain  
$\mu_3 \bigl(\ST_N \bigr) = 
(\mu_2 \circ \rho )\bigl(\ST_N \bigr) = \ST_N$.

Thus, we have obtained all relations given in
Theorem~\ref{main}.

%%%%%%%%%%%%%%%%%%%%%%%%%%%%%%%%%%%%%%%%%%%%%%%%%%%%%%%%%%%%%%%%%%%%%%%%
\appendix
\section{Appendix}
\subsection{Proofs of propositions of Section~2.1}

{\bf Proof of Lemma~\ref{defWalt}.}\ 
Let $\bigl\lfloor W \bigr\rfloor_{N+1} \in {\mathcal B}(n,N+1)$ 
denote the word which is
obtained from a word $W \in {\mathcal B}(n-1,N)$ 
by the replacement
$R_{a_1,\ldots,a_{n-1}} \to R_{a_1,\ldots,a_{n-1},{N+1}}$
applied to all factors. A key step in proving
Lemma~\ref{defWalt} is to observe that
the following recursive relation holds
\begin{equation}\label{Wp}
     W(n,N) \, 
    \bigl\lfloor W(n-1,N) \bigr\rfloor_{N+1}  
    = W(n,N+1) \,.
\end{equation}
For instance,
\begin{align}\nonumber
{}&  
 W(2,3) \bigl\lfloor W(1,3) \bigr\rfloor_{4} =
 R_{12} R_{13} R_{23} \bigl\lfloor R_{1} R_{2} R_3 \bigr\rfloor_4  = 
 R_{12} R_{13} R_{14} R_{23}  R_{24} R_{34} = W(2,4) .
\end{align}  

For $n=2$, eq. (\ref{Wp}) is almost trivial since 
$\lfloor W(1,N) \rfloor_{N+1} = R_{1,N+1} \ldots R_{N,N+1}$
and for each factor $R_{a,N+1}$, $a\neq N$  the rightmost
factor in $W(2,N)$ with which it does not commute
is $R_{a,N}$. Therefore, moving all factors from 
$\lfloor W(1,N) \rfloor_{N+1}$ to the left until
they meet their non--commuting counterparts, we
obtain the lexicographically ordered word $W(2,N+1)$.
Note that only the last two factors in 
$\lfloor W(1,N) \rfloor_{N+1}$ have the same 
non--commuting counterpart, namely, $R_{N-1,N}$.

The proof for $n \geq 3$ is similar.
Each factor $R_{a_1,\ldots,a_{n-1},N+1}$
in $\lfloor W(n-1,N) \rfloor_{N+1}$ belongs
to a cluster, i.e. a product of consecutive 
lexicographically ordered factors. 
A cluster has length one if $a_{n-1} < N-1$.
If $a_{n-1} \geq N-1$, then the length of
a cluster is $k+1$, where $k=\max\{m :  a_{n-m} \geq N-m \}$.
If $R_{a_1,\ldots,a_{n-k-1},N-k,\ldots,N-1,N+1}$ is the
leftmost factor of a cluster of length $k+1$, 
then the rightmost non--commuting counterpart in $W(n,N)$ 
for all factors of the cluster is 
$R_{a_1,\ldots,a_{n-k-1},N-k,\ldots,N-1,N}$ 
Therefore, moving each cluster from 
$\lfloor W(n-1,N) \rfloor_{N+1}$
to the left until it meets its counterpart, we achieve
the lexicographic order of all factors, that is,
we obtain the word~\hbox{$W(n,N+1)$}.

The recursive relation (\ref{Wp}) implies that
\begin{equation}\label{Wpp}
 W(n,N) = 
 \prod_{n\leq k \leq N}^{\rightarrow} 
 \bigl\lfloor W(n-1,k-1)\bigr\rfloor_{k} \,,
\end{equation}
which in turn leads to the   
expression (\ref{Wexp}) for~$W'(n,N)$.
\qed\vspace{1mm}

To prove Theorem~\ref{multitetra}
we will need also the following statement.
\begin{lem}\label{LEX} 
Let $\omega$ be an involutive automorphism of
${\mathcal B}(n,N)$ such that 
\begin{equation}\label{om}  
  \omega \bigl( R_{a_1,\ldots,a_{n}} \bigr) =
  R_{N+1-a_{n},\ldots,N+1-a_{1}}  \,, 
\end{equation}
for all generators of~${\mathcal B}(n,N)$. 
Then for all $N \geq n$ the equalities
\begin{equation}\label{omWW}  
  \omega \bigl( W(n,N) \bigr)  = 
   \eta \bigl( W'(n,N) \bigr) \,, \qquad
   \omega \bigl( W'(n,N) \bigr)  = 
   \eta \bigl( W(n,N) \bigr) 
\end{equation}
hold in the strong sense, that is their r.h.s.
coincide with their l.h.s. without a permutation 
of commuting factors.
\end{lem}

\noindent {\bf Proof of Lemma~\ref{LEX}.}\
Denote $\bar{W}'(n,N)=\eta \bigl( W'(n,N)\bigr)$.
Applying the anti--automorphism
$\eta$ to (\ref{Wexp}) and (\ref{Wp}), we obtain  
\begin{align}
\label{barWexp}  
{}&  \bar{W}'(n,N) = \!
 \mathop{{\prod}'}\limits^{\leftarrow}_{1\leq b_1< \ldots < b_n \leq N }  
 \! R_{b_1,\ldots,b_{n}}    \,,\\
\label{Wm}
{}&  \bar{W}'(n,N+1) = 
   \bigl\lfloor \bar{W}'(n-1,N) 
  \bigr\rfloor_{N+1} \, \bar{W}'(n,N) \,.
\end{align} 
Applying to (\ref{barWexp}) the automorphism
$\omega$, and relabelling  the indices by
$a_k = N+1 - b_{n+1-k}$, we recover the r.h.s of 
formula (\ref{lex}). Hence we infer that
$\omega \bigl( \eta \bigl( W'(n,N) \bigr) \bigr) = W(n,N)$,
which proves the first equality in~(\ref{omWW}).
The second equality follows then immediately
since $\omega$ and $\eta$ mutually
commute and both are involutive.
\qed\vspace{2mm}

\noindent {\bf Proof of Theorem~\ref{multitetra}.}\
For $n=1$, the statement of the Theorem is obvious
since all $R$'s commute.
For $n=2$, the statement can be proven by 
induction on~$N$. The base, for $N=3$, is simply the
relation (\ref{Rgen}).
Assume that the equality 
$W(2,N)=\bar{W}(2,N)$ has been already established
for some $N \geq 3$. Then, taking into account the
recursive structure of (\ref{Wp}) and (\ref{Wm}), 
we have to prove that
\begin{equation}\label{W22} 
 \bigl\lfloor \prod_{1 \leq a \leq N }^{\rightarrow} 
 R_{a} \bigl\rfloor_{N+1} = 
  \bigl( W(2,N) \bigr)^{-1} 
  \bigl\lfloor \prod_{1 \leq a \leq N }^{\leftarrow} 
  R_{a} \bigl\rfloor_{N+1} 
  \,  W(2,N)  \,. 
\end{equation}
Observe that the Yang--Baxter equation (\ref{RYB}) can be 
rewritten as the following ``almost commutativity" relation
for two $R$'s,
\begin{equation}\label{RYB2}  
  \bigl\lfloor R_a \, R_b  \bigl\rfloor_{N+1} = 
  R^{-1}_{ab} \,\bigl\lfloor R_{b} \, R_{a} \bigl\rfloor_{N+1} R_{ab} \,.
\end{equation}
Using this relation, we can move $R_{1,N+1}$ in 
the   product 
$\prod\limits^{\rightarrow} R_{a, N+1} $ to the right,
then move $R_{2,N+1}$, and so on until we obtain 
the reverse  ordered product
$\prod\limits^{\leftarrow} R_{a, N+1} $.
Notice that all the extra ``twisting" factors 
$R^{\pm 1}_{ab}$ arising
in this process commute with any $R_{c,N+1}$,
$c\neq a,b$. 
Therefore, at each step these twisting factors can be moved
outside of the product 
$\bigl\lfloor \ldots \bigl\rfloor_{N+1}$. It is easy
to see that these factors combine into 
$\bigl(W(2,N)\bigr)^{\pm 1}$. 
Thus, we have established the inductive step~(\ref{W22})
and hence the Theorem is proven for~$n=2$.

For $n=3$, the Theorem can be proven along the same lines.
The base, for $N=4$, is the relation (\ref{Rgen}).
Assuming that 
$W(3,N)=\bar{W}(3,N)$ has been already established
for some $N \geq 4$ and taking into account 
(\ref{Wp}) and (\ref{Wm}), we have to prove the inductive step
\begin{equation}\label{W33}  
  \bigl\lfloor W(2,N) \bigl\rfloor_{N+1} = 
  \bigl( W(3,N) \bigr)^{-1} 
  \bigl\lfloor \bar{W}(2,N) \bigl\rfloor_{N+1} \,
  W(3,N) \,.
\end{equation}
Observe that the tetrahedron equation (\ref{RTE}) can
be rewritten as the following ``almost Yang--Baxter" relation
for three $R$'s,
\begin{equation}\label{RTE22}  
  \bigl\lfloor R_{ab} \, R_{ac} \, R_{bc}  \bigl\rfloor_{N+1} =
  R_{abc}^{-1} \, 
  \bigl\lfloor R_{bc} \, R_{ac} \, R_{ab}  \bigl\rfloor_{N+1} \, R_{abc} \,,
\end{equation}
where $a<b<c$. Notice that $R_{abc}$ commutes
with any factor in $\bigl\lfloor W(2,N) \bigl\rfloor_{N+1}$ 
except those on the l.h.s of (\ref{RTE22}). 
Therefore, the factors in $\bigl\lfloor W(2,N) \bigl\rfloor_{N+1}$
can be reordered exactly in the same way as the factors
in $W(2,N)$ and the extra twisting factors $R_{abc}^{\pm 1}$
arising at each step can be moved 
outside of the product $\bigl\lfloor \ldots \bigl\rfloor_{N+1}$.
It is easy to see that these factors combine into 
$\bigl(W(3,N)\bigr)^{\pm 1}$. Since we already know that
$W(2,N)=\bar{W}(2,N)$, it follows that the inductive
step (\ref{W33}) holds
and so the Theorem is proved for~$n=3$.

It is now clear that the proof continues by the double
induction on $n$ and~$N$. 
For a given $n$, relation (\ref{Rgen}) can be rewritten
as ``almost" the relation for $n-1$. Using it,
we can reorder the factors in 
$\bigl\lfloor W(n-1,N) \bigl\rfloor_{N+1}$ exactly in the same
way as the factors in $ W(n-1,N)$.
At each step the extra twisting factors 
can be moved outside of the product 
$\bigl\lfloor \ldots  \bigl\rfloor_{N+1}$ due to relation~(\ref{Rcomm}).
Therefore, the same line of arguments as for $n=2,3$
proves the Theorem for an arbitrary~$n$.
\qed\vspace{2mm}

\noindent {\bf Proof of Lemma~\ref{TW23}.}\\
$\bullet$\
The first equality in (\ref{TWRZa}) is obvious from
Definition~\ref{defT}  and eq.~(\ref{Rkdef}). 
In order to prove the second equality in (\ref{TWRZa}), 
we consider the following elements of 
${\mathcal T}_{N}$
\begin{equation}\label{Trm}
  \ST_{N,b} =   
\prod_{1\leq a < k \leq b }^{\rightarrow}   
  \ddil{ ak}   \,, \qquad  2 \leq b \leq N .
\end{equation}
We observe that 
\begin{equation}\label{Trb}
  \ST_{N,b} =    \ST_{N,b-1}
  \prod_{1 \leq a < b}^{\rightarrow} \ddil{ ab } \,.
\end{equation}
Indeed, it follows from the definition (\ref{Rkdef})
that the rightmost factor in the lexicographically
ordered product $\ST_{N,b}$ which does not commute with
$\ddil{ ab}$ is $\ddil{ a,b-1}$. 
Repeating the argument used in the proof of
formula (\ref{Wp}) for $n=2$, we conclude that
(\ref{Trb}) holds. And it is easy to see that (\ref{Trb})
along with 
$\ST_N = \ST_{N,N}$ implies that 
\begin{equation}\label{Trbb}
  \ST_N = \ST_{N,N-1}
  \prod_{1 \leq a < N}^{\rightarrow} \ddil{ aN } 
  = \ldots =
  \mathop{{\prod}'}\limits^{\rightarrow}_{1\leq a< b \leq N }    
  \ddil{ab}
   = \bigl( W'(2,N) \bigr)^\sharp \,.
\end{equation} 

\noindent $\bullet$\
The first equality in (\ref{TWRZb}) is obvious from
Definition~\ref{defT}.
The second equality in (\ref{TWRZb}) is equivalent 
to the statement that  
\begin{equation}\label{TrR}
  \ST_N = \ST_{N-1} \,  \SW'(2,N) \,,
\end{equation} 
where $ \SW'(2,N)$ is given by 
(\ref{W2lex}) with each $R_{ab}$ replaced with 
$\ddil{a,b,N+1}$.
%% $\dil{ \SZ_{ab }  }^{(N-b)}$.
Using relations (\ref{Wp}) and (\ref{Trbb}),  
we can rewrite the r.h.s. of (\ref{TrR}) as follows
\begin{equation}\label{TrRR}
  \tilde{\SW}'(2,N-1) \, \SW'(2,N-1)  \, \Bigl(
  \prod_{1 \leq e < N}^{\rightarrow}  \ddil{ eN } \Bigr) \,,
\end{equation} 
where $ \tilde{\SW}'(2,N)$ is given by 
(\ref{W2lex}) with each $R_{ab}$ replaced with  
$\ddil{ ab }'$ which is an element of ${\mathcal T}_{N-1}$ 
given by~(\ref{Rkdef}). 
Definition (\ref{Rkdef}) implies that $\ddil{ ab }'$ 
is the rightmost
factor in $\tilde{\SW}'(2,N-1)$ which does not commute
with the factor $\ddil{a,b,N+1}$ in 
$\SW'(2,N-1)$. Therefore, moving all factors from
$\SW'(2,N-1)$ to the left until
they meet their non--commuting counterparts
and noticing that 
$\ddil{ ab}' \ddil{a,b,N+1}=
\ddil{ ab }$, we conclude that
(\ref{TrRR}) coincides with the second expression
in (\ref{Trbb}) and hence relation (\ref{TrR}) 
holds. 

\noindent $\bullet$\
To prove the last part of the Lemma we
consider $ \SW'(2,N)$ introduced in (\ref{TrR}). 
One can check that all its factors containing
$\SZ_{N-1,N}$ can be moved to the right preserving their 
order. Making then the same procedure with
factors containing $\SZ_{N-2,N}$, $\SZ_{N-3,N}$, etc.,
we conclude that $ \SW'(2,N)$ is given by the same expression 
in (\ref{W2lex}) where each $R_{ab}$ is replaced with 
$\ddil{a,N+1+a-b,N+1}$.
%% $\dil{ \SZ_{a,N+1+a-b }  }^{(b-a-1)}$.
With such a form of $ \SW'(2,N)$ formula (\ref{TrR})
leads immediately to the second equality in~(\ref{TWRZc}).

Equivalence of the first and second expressions for $\ST_N$
in (\ref{TWRZc}) is verified as follows
\begin{equation}\label{wwmm}
  \bigl( W(3,N)\bigr)^{\star\star}
  = \mu_2 \Bigl( \bigl( W'(3,N)\bigr)^{\star} \Bigr)
  = \mu_2 \Bigl( \bigl( W(3,N)\bigr)^{\star} \Bigr)
  = \bigl( W'(3,N)\bigr)^{\star\star} .
\end{equation} 
The middle equality is due to the second part of the Lemma.
In the first and last equalities we used Lemma~\ref{mu2T}, 
which is proven below in Section~A.4. 
\qed

\subsection{Proofs of propositions of Section~2.2}
 
\noindent {\bf Proof of Theorem~\ref{center}.}\
Recall that 
 $\alpha \in {\mathbb Z}_{\geq 0}^{N(N-1)/2}$
is a vector comprised of the weights $\alpha_{ij}$ 
assigned to the vertices of $Q_N$ ordered lexicographically.
By (\ref{mona}), each $\alpha$ determines a 
monomial $\SM(\alpha) \in {\mathcal T}_N$.
Relations (\ref{ZZ}) imply that
$\SM(\alpha)$ commutes with
all generators of ${\mathcal T}_N$ iff 
$B \alpha =0$.
Thus, we have to study the 
kernel of the incidence matrix~$B$.  

$\bullet$ 
Let $\mathcal U$ be the symmetry transformation of the vertices 
of $Q_N$ such that ${\mathcal U} (i,j)=(j-i,N+1-i)$, and 
$U$ be the matrix of the corresponding 
orthogonal transformation of the basis in 
${\mathbb Z}_{\geq 0}^{N(N-1)/2}$. 
Applying $\rho$ to (\ref{ZZ}), we infer that
$B$ commutes with $U$. Taking
into account that $U$ is orthogonal and 
$B$ is skew--symmetric, we conclude that if 
$\alpha \in \mathrm{Ker} B$,
then $\alpha$ is an eigenvector of~$U$.
However, $U^3=1$ and so
the only real eigenvalue of $U$ is one.
Thus, $\SM(\alpha) \in {\mathcal Z}_N$ implies that
$U \alpha = \alpha$. Hence 
\begin{equation}\label{alrho}
 \alpha_{{\mathcal U} (i,j)} = \alpha_{ij}  \,,
\end{equation} 
which is the first symmetry in (\ref{alsym}).
Applying $\rho$ to a central monomial
$\SM(\alpha)$ and taking this symmetry of
its weights into account, we infer that
$\rho \bigl( \SM(\alpha)\bigr) = 
q^{\varepsilon(\alpha)}  \SM(\alpha)$,
where the multiplicative constant appears due
to reordering of the generators. However,
the property $\rho \circ \rho \circ \rho =id$
implies that $q^{\varepsilon(\alpha)}=1$ and
thus $\rho \bigl( \SM(\alpha)\bigr) = \SM(\alpha)$.

Let $\mu$ be any of the anti--homomorphisms in 
(\ref{mudef}), $\mathcal K$ be the corresponding 
symmetry transformation of the vertices of $Q_N$
(e.g. ${\mathcal K}(i,j)=(j-i,j)$ for $\mu_1$), and 
$K$ be the matrix of the corresponding 
orthogonal transformation of the basis in 
${\mathbb Z}_{\geq 0}^{N(N-1)/2}$.
Applying $\mu$ to (\ref{ZZ}), we infer that 
$K B K = - B$.
Hence $B^2$ commutes with $K$.
Taking into account that both these matrices are
real symmetric and $K^2=1$,  
we conclude that if 
$\alpha \in \mathrm{Ker} B$,
then $K \alpha = \pm \alpha$.
But $\alpha_{12}=\alpha_{1N}=\alpha_{N-1,N}$
thanks to~(\ref{alrho}).
Therefore, $K \alpha = \alpha$. Hence
\begin{equation}\label{alrefl}
 \alpha_{{\mathcal K}(i,j)} = \alpha_{ij}  \,,
\end{equation}
which yields the remaining symmetries in~(\ref{alsym}).
Applying $\mu$ to a central monomial
$\SM(\alpha)$ and taking these symmetries 
into account, we obtain that
$\mu \bigl( \SM(\alpha)\bigr) = 
q^{\varepsilon'(\alpha)}  \SM(\alpha)$.
The property $\mu \circ \mu =id$
implies that $q^{\varepsilon'(\alpha)}=1$ and
thus $\mu \bigl( \SM(\alpha)\bigr) = \SM(\alpha)$.

$\bullet$
Consider monomials $\SM(\alpha)$
and $\SM(\alpha')$ such that 
$\alpha_{1j} = \alpha'_{1j}$ for $j=2,\ldots,N$. 
Then, in view of (\ref{alrho}), all components of
vector $\alpha''=\alpha-\alpha'$ corresponding
to boundary vertices of $Q_N$ vanish. 
Suppose that $\SM(\alpha)$
and $\SM(\alpha')$ are central. Then
$\SM(\alpha'')$ commutes with all generators 
of~${\mathcal T}_N$. In particular, the condition
that $\SM(\alpha'')$ commutes
with all $\SZ_{1j}$ is equivalent to
a system of equations:
 $\alpha''_{2,j+1}=\alpha''_{2,j}$.
However, we already have $\alpha''_{23}=0$. 
Therefore, we conclude that $\alpha''_{2j}=0$
for $j=3,\ldots,N$. This implies, in view of (\ref{alrho}),
that all components of $\alpha''$ corresponding
to next to boundary vertices of $Q_N$ also vanish. 
Continuing this consideration similarly
for $\SZ_{2j}$, $\SZ_{3j}$ etc., we conclude
that $\alpha''=0$. Thus, $\SM(\alpha')=\SM(\alpha)$,
i.e. two central monomials coincide iff
they have coinciding weights at one boundary of~$Q_N$.
Thanks to the symmetry (\ref{alrefl}), the 
latter condition is equivalent to a weaker 
condition:  $\alpha'_{1j} = \alpha_{1j}$ 
for $j=2,\ldots,\chi(N)+1$.

$\bullet$
Let us show that, given arbitrary integers 
$\beta_1, \ldots, \beta_{\chi(N)}$, there exists
a unique integer vector $\tilde{\alpha}$ such
that $\tilde{\alpha}_{1j}=\beta_{j-1}$ for 
$j=2,\ldots,\chi(N)+1$
and $\SM(\tilde{\alpha})$ given by (\ref{mona}) 
(where some $\tilde{\alpha}_{ij}$ can be negative)
commutes with all generators $\SZ_{ij}$.
Let $\partial Q_N$ stand for the set of all
 boundary vertices of $Q_N$. Forgetting about
the edges, we have $Q_N = 
\partial Q_N \cup  \partial Q_{N-3} \cup \ldots$.

First, given the weights $\tilde{\alpha}_{1j}$
for $j=2,\ldots,\chi(N)+1$, we extend them 
to weights at other vertices of $\partial Q_N$
by the symmetries (\ref{alrho}) and (\ref{alrefl}).
Now, the requirement that $\SM(\tilde{\alpha})$
commutes with all generators assigned to $\partial Q_N$
 fixes uniquely all weights 
at $\partial Q_{N-3}$. Indeed, $\SM(\tilde{\alpha})$
commutes with $\SZ_{1j}$ if $\tilde{\alpha}_{2,j+1}-
\tilde{\alpha}_{2j}=\tilde{\alpha}_{1,j+1}-
\tilde{\alpha}_{1,j-1}$. 
  Taken into 
account that $\tilde{\alpha}_{23}=\tilde{\alpha}_{13}$,  
this set of equations determines all $\tilde{\alpha}_{2j}$
uniquely. Moreover, the symmetry (\ref{alrefl}) for 
$\tilde{\alpha}_{1j}$ induces the same symmetry
for $\tilde{\alpha}_{2j}$. Finally, $\tilde{\alpha}_{2j}$
can be extended to weights at other vertices of 
$\partial Q_{N-3}$ by the symmetry~(\ref{alrho}).

Similarly, given weights at $\partial Q_{k+3}$
and $\partial Q_{k}$, the requirement that $\SM(\tilde{\alpha})$
commutes with all generators assigned to $\partial Q_k$
 fixes uniquely all weights at $\partial Q_{k-3}$. 
 Indeed, $\SM(\tilde{\alpha})$
commutes with all $\SZ_{ij}$ for a given $i$ if 
$\tilde{\alpha}_{i+1,j+1}-
\tilde{\alpha}_{i+1,j} =\tilde{\alpha}_{i,j+1}
-\tilde{\alpha}_{i,j-1}+\tilde{\alpha}_{i-1,j-1}-
\tilde{\alpha}_{i-1,j}$. These equations
are resolved uniquely since we know the r.h.s.
and several first values of $\tilde{\alpha}_{i+1,j}$
thanks to the symmetry~(\ref{alrho}).

Thus, the vector $\tilde{\alpha}$ is recovered
uniquely {}from its first $\chi(N)$ components. However,
it can happen that some $\tilde{\alpha}_{ij}$ are negative.
In this case we take another vector,
$\alpha=\tilde{\alpha} +m\, \alpha^0$, 
where $\alpha^0$ is a vector such that $\alpha^0_{ij}=1$,
and $m$ is sufficiently large positive number
to ensure positivity of all  $\alpha_{ij}$. 
Then we have
$\SM(\alpha) = q^{\varepsilon} \cdot
 \SM(\tilde{\alpha}) \bigl( \SM(\alpha^0) \bigr)^m$.
It remains to observe that $\SM(\alpha^0) \in {\mathcal Z}_N$
and hence $\SM(\alpha) \in {\mathcal Z}_N$.
\qed\vspace{2mm}

\subsection{Proofs of propositions of Section~2.3}
\noindent {\bf Proof of Lemma~\ref{RSR}.}\
The statement of the Lemma for $\SF$ follows from
that the action of
$\SF_{123} \, \SF_{124} \, \SF_{134} \, \SF_{234}$ and
$\SF_{234} \, \SF_{134} \, \SF_{124} \, \SF_{123}$
on the monomial 
$x_{12}^a x_{13}^b x_{23}^c x_{14}^d x_{24}^e x_{34}^f$,
yields the same result.

Note that $\SR$ given by (\ref{Rdef}) is invertible.
If the argument of the quantum exponential in 
(\ref{Rdef}) is denoted by $\SX$, then we have
$\SR^{-1} = \dil{\SX}^{-1} \SF$, where 
$ \dil{x}^{-1}$ is the following formal series
\begin{equation}\label{Sxinv}
  \dil{x}^{-1} =  \sum_{n \geq 0}   
   \frac{ q^{n(n-1)/2}\, x^n}{(1-q)\ldots(1-q^{n})} \,.
\end{equation}

In order to check that $\SR$ satisfies the
tetrahedron equation, one has to substitute $\SR_{abc}$
into (\ref{RRR}), move all
$\SF$'s to the left, and cancel the products of 
$\SF$'s on the both sides of the equation by invoking
the first part of the Lemma. Then one is left 
with the equality
\begin{equation}\label{4444}
\dil{\SX} \, \dil{\SX\SY} \, \dil{\SZ} \, \dil{\SY}  =  
\dil{\SZ} \, \dil{\SZ\SX} \, \dil{\SY} \, \dil{\SX}  \,,
\end{equation} 
where 
\begin{equation}\label{xyz44}
\begin{aligned}
{}& \SX = q^{1+\gamma} \sx_{12} \sx_{23}^{-1}
 \sy_{12}^{\gamma} \sy_{13}^{-\gamma} \sy_{23}^{-1-\gamma} , \quad
\SY= q^{1+\gamma} \sx_{23} \sx_{34}^{-1}
\sy_{23}^{\gamma} \sy_{24}^{-\gamma} \sy_{34}^{-1-\gamma} , \\[1mm]
{}& \qquad\qquad\qquad  
 \SZ = q^{1+\gamma}    \sx_{13} \sx_{24}^{-1} 
\sy_{13}^{\gamma} \sy_{14}^{-\gamma} 
\sy_{23}^{1+\gamma} \sy_{24}^{-1-\gamma} . 
\end{aligned}
\end{equation}
These operators satisfy relations~(\ref{XYZ}).
Therefore, comparing the first equality in (\ref{44}) 
with (\ref{4444}), we conclude that equality 
(\ref{4444}) holds. 
\qed\vspace{2mm}

\noindent {\bf Proof of Lemma~\ref{Zrep}.}\   
Note that 
$\tau(\SZ_{ab})$ does not commute with $\tau(\SZ_{a'b'})$ 
only in the following six cases: 1) $a'=a$ and 
$b'=b \pm 1$; 2) $a'-a=\pm 1$ and
$b'=b$ or $b'=b+a'-a$.
It is easy to check that in these cases 
$\tau(\SZ_{ab})$ and $\tau(\SZ_{a'b'})$ satisfy the defining
relations~(\ref{ZZ}). 

Clearly, the representation (\ref{Zxy}) is faithful
for ${\mathcal T}_{2}$ since in this case the algebra
has only one generator.
The proof of faithfulness of this representation in
other cases will use the inclusion
${\mathcal T}_{N-1} \subset {\mathcal T}_{N}$
(cf.~Remark~\ref{AAA}). 
Let $N$ be the minimal positive number such that
the representation $\tau$ is not faithful for ${\mathcal T}_{N}$.
Then there exists a polynomial $H$ in ${N} \choose {2}$ 
variables such that 
$\SH \equiv H(\tau(\SZ_{12}),\ldots,\tau(\SZ_{N-1,N}))=0$.

Without a loss of generality we can assume that
$\SH$ has the form
\begin{equation}\label{PZ}    
 \SH =
 \sum_{k_1, \ldots, k_{N-1} \geq 0} B_{k_1 \ldots k_{N-1}} \,
 \tau\bigl(\SZ_{1N}^{k_1} \ldots \SZ_{N-1,N}^{k_{N-1}} \bigr) \,
 \SH_{k_1 \ldots k_{N-1}}  \,,
\end{equation}
where $B_{k_1 \ldots k_{N-1}} \neq 0$ at least
for one set $\{k_1, \ldots, k_{N-1}\}$ such that
$k_1 + \ldots + k_{N-1} >0$.
$\SH_{k_1 \ldots k_{N-1}}$ stand for  
polynomials in the generators of~${\mathcal T}_{N-1}$ evaluated
in the representation~(\ref{Zxy}).

Acting by $\SH$ on $f \in S_{N} \subset S_{N+1}$,
we obtain
\begin{equation}\label{QZ}    
 \SH \cdot f =
 \sum_{k_1 \ldots k_{N-1} \geq 0}   
 q^{\varepsilon(k_1,\ldots,k_{N-1})} B_{k_1 \ldots k_{N-1}} \,
 \biggl(
 \prod_{a=1}^{N-1} \Bigl( \frac{x_{aN}}{x_{a+1,N+1}} \Bigr)^{k_a} 
 \biggr) \,
 \Bigl( \tilde{\SH}_{k_1 \ldots k_{N-1}} \cdot f \Bigr) \,,
\end{equation}
where $\tilde{\SH}_{k_1 \ldots k_{N-1}}
 =  \SY_N \SH_{k_1 \ldots k_{N-1}}$ with
$\SY_N=\sy_{2N} \ldots \sy_{N-1,N}$.
Note that all terms in the sum in (\ref{QZ}) are 
linearly independent monomials in 
$x^{-1}_{2,N+1},\ldots,x^{-1}_{N,N+1}$.
Hence, $\SH=0$ implies
that $\tilde{\SH}_{k_1 \ldots k_{N-1}} \cdot f =0$.
Since $\SY_N$ is invertible, we conclude that 
$\SH_{k_1 \ldots k_{N-1}}$ annihilates an
arbitrary $f$ from~$S_N$ and, thus,
$\SH_{k_1 \ldots k_{N-1}}=0$. But this
contradicts the assumption that the representation
$\tau$ is faithful for ${\mathcal T}_{N-1}$.
\qed\vspace{2mm}

\subsection{Proof of propositions of Section~2.4 }

\begin{lem}\label{mu2T} 
For all $N \geq 3$, the following equalities
\begin{equation}\label{muW}
 \mu_2 \Bigl( \bigl( W(3,N)\bigr)^{\star} \Bigr)
  = \bigl( W'(3,N)\bigr)^{\star\star} , \qquad
  \mu_2 \Bigl( \bigl( W'(3,N)\bigr)^{\star} \Bigr)
  = \bigl( W(3,N)\bigr)^{\star\star} ,
\end{equation} 
hold in the strong sense that is their r.h.s. coincide
with their l.h.s. without a permutation of commuting
factors.
\end{lem} 

\noindent {\bf Proof of Lemma~\ref{mu2T}.}\
The first equality in (\ref{muW})
is checked with the help of Lemma~\ref{LEX} as follows:
\begin{align*}
{}& \mu_2 \Bigl( \bigl( W(3,N) \bigr)^\star \Bigr) = 
 \mu_2 \Bigl(
 \prod_{1\leq a < b < c \leq N}^{\rightarrow}  
  (R_{abc})^\star \Bigr) 
  =
 \mathop{{\prod}'}\limits^{\rightarrow}_{1\leq a< b <c \leq N }   
  \mu_2 \Bigl( \bigl(\omega (R_{abc}) \bigr)^\star \Bigr) \\
{}& =  \mathop{{\prod}'}\limits^{\rightarrow}_{1\leq a< b <c \leq N } 
  \mu_2 \Bigl( \ddil{ N+1-c,N+1-b,N+1-a }   \Bigr) 
 =    \mathop{{\prod}'}\limits^{\rightarrow}_{1\leq a< b <c \leq N } 
 \ddil{ a,c+a-b,c } 
 =  \bigl( W'(3,N) \bigr)^{\star\star} .
\end{align*}  
The second equality can be checked in a similar way.
\qed\vspace{2mm}

Recall that the homomorphisms 
$\theta$ and $\phi_\gamma$ were introduced 
in Lemma~\ref{RSR}.
In this subsection we will use the following notations: 
\begin{equation}\label{Rabc}    
 \SR_{abc} = \phi_0(R_{abc}) = \SF_{abc} \cdot 
 \check{\SR}_{abc} , \qquad
 \SF_{abc} = \theta \bigl( R_{abc} \bigr) , \qquad
 \check{\SR}_{abc} =
 \dil{ q \sx_{ab}   \sx_{bc}^{-1} 
    \sy_{bc}^{-1} }   .
\end{equation}
We will need also the following homomorphisms 
$\varphi, \varphi'$ from
${\mathcal B}(3,N)$ to $\text{End}\,\bigl(S_N\bigr)$:
\begin{align}
\label{phi}  
{}& \varphi ( R_{abc} ) =
 \dil{ \frac{q \sx_{ab} }{ \sx_{a+c-b,c} }
 \prod_{k=1}^{c-b} \frac{ \sy_{a+k,b+k-1}}{\sy_{a+k,b+k}}   }  
 , \\
\label{phihat} 
{}&  \varphi' ( R_{abc} )  =  
 \dil{ q \Bigl(
 \prod_{k=1}^{c-b} \frac{\sx_{b-a,b+k-1} }{ \sx_{b-a+1,b+k} }
 \Bigr)
   \frac{\sy_{b-a+1,b}}{\sy_{b-a+1,c}} } ,
\end{align} 
where $\sy_{ab} \equiv 1$ if $a=b$.

Recall that the word $W'(3,N)$ is given by (\ref{W3lex}),
and that $\bar{W}'(3,N)=\eta \bigl(W'(3,N)\bigr)$
is the word with reversely ordered factors.
Let us introduce
\begin{equation}\label{Fr}  
 \SF_N = \theta \bigl(W'(3,N)\bigr) = 
  \theta \bigl(\bar{W}'(3,N)\bigr) ,
\end{equation}
where the equality of the expressions is due to
Theorem~\ref{multitetra}.

\begin{lem}\label{RR'}
For all $N \geq 3$, we have the equalities
\begin{align}
\label{RFF}
{}&  \phi_0 \bigl(W'(3,N)\bigr) = 
  \SF_N \cdot
  \varphi \bigl(W'(3,N)\bigr) \,, \\[1mm]
\label{RFFr} 
{}&   \phi_0 \bigl(\bar{W}'(3,N)\bigr) = 
  \SF_N \cdot
   \varphi' \bigl(\bar{W}'(3,N)\bigr) \,.
\end{align} 
\end{lem}

\noindent {\bf Proof of Lemma~\ref{RR'}.}\
Eq. (\ref{Wpp}) implies that
\begin{equation}\label{w3t}  
 \phi_0 \bigl( W'(3,N) \bigr) = 
 \prod_{3\leq c \leq N}^{\rightarrow} 
 \phi_0 \Bigl( \bigl\lfloor W'(2,c-1)\bigr\rfloor_{c} \Bigr) \,.
%% = \phi_0 \bigl( W'(3,N-1) \bigr) \,
%% \phi_0 \Bigl( \bigl\lfloor W'(2,N-1)\bigr\rfloor_{r} \Bigr) .
\end{equation}
Consider $\check{\SR}_{abc}$ entering 
the $\phi_0\bigl(\lfloor\ldots\rfloor_c \bigr)$ part.
In order to establish (\ref{RFF}) it 
suffices to prove that:\\
{\bf i}) $\check{\SR}_{abc}$ 
transforms into $\varphi(\SR_{abc})$ when all $\SF$'s
from $\phi_0\bigl(\lfloor\ldots\rfloor_c \bigr)$ 
are moved to the left of 
this part;\\ 
{\bf ii})  $\varphi(\SR_{abc})$ does not change
when, for any $k>c$, all $\SF$'s from 
$\phi_0\bigl(\lfloor W'(2,k-1)\rfloor_k \bigr)$
are moved through it.

{\bf Proof of i).}\
Consider $\check{\SR}_{abc}$ which is not
the rightmost one (otherwise, it already
has the form $\varphi(R)$). 
All $\SF$'s to the right of it have the form 
$\SF_{a'b'c}$, $\{a',b'\}\neq \{a,b\}$. 
By (\ref{Fxx}), they do not change $\sx_{ab}$
but act on the first index of $\sx_{bc}$ 
as permutations $P_{a'b'}$. {}From the recursive
structure of the word 
$W'(2,c-1)=\lfloor W(1,1)  \rfloor_{2} \ldots
  \lfloor W(1,c-2)  \rfloor_{c-1} $ it
follows that $\sx_{bc}$ has to
be pulled through the following chain of permutations: 
$$
P_{a+1,b}\ldots P_{b-1,b} \ P_{1,b+1} \ldots P_{b,b+1} \ldots
 P_{1,c-1} \ldots P_{c-2,c-1} . 
$$
The part  $P_{a+1,b}\ldots P_{b-1,b}$ here is present
if $a<b-1$. It transforms
$\sx_{bc}$ into $\sx_{a+1,c}$. The remaining
permutations are grouped into $(c-1-b)$ shift operators
$U_d =P_{d,d+1} \ldots P_{12}$,\  $b \leq d \leq c-2$. 
Each $U_d$ increases the first index of 
$\sx_{a+1,c}$ by one, thus, transforming it into
$\sx_{a+c-b,c}$.

Consider now $\sy_{bc}$ entering the argument 
of~${\check\SR}_{abc}$. It is first transformed 
by $\SF_{a+1,b,c}$ into $\frac{\sy_{a+1,c}}{\sy_{a+1,b}}$.
Note that, by (\ref{Fyy}), all $\SF$'s entering
$\phi_0\bigl(\lfloor W'(2,c-1) \rfloor_{c}\bigr)$ 
act non--trivially
only on $\sy$'s whose second index is~$c$.
{}From the recursive structure of the word $W'(2,c-1)$ it
follows that $\sy_{a+1,c}$ is transformed by consecutive
action of pairs $\SF_{a+k,b+k,c} \SF_{a+1+k,b+k,c}$,
$1 \leq k \leq c-b-1$. Each such pair
transforms $\sy_{a+k,c}$ into 
$\frac{\sy_{a+k,b+k}}{\sy_{a+1+k,b+k}} \sy_{a+1+k,c}$,
which yields the product of $\sy$'s in~(\ref{phi}). 

{\bf Proof of ii).}\ 
This part is trivial for $\sy$'s since, by 
(\ref{Fyy}), $\sy_{ab}$ commutes with any
$\SF_{a^{\prime} b^{\prime} k}$ if $k \neq a,b$.
For $\sx$'s we have to consider the transformations of
$\frac{\sx_{ab} }{ \sx_{a'c} }$, where $a' \equiv a+c-b$.
Note that $a< a' <c$. Therefore, 
$\frac{\sx_{ab} }{ \sx_{a'c} }$ is first transformed
by $\SF_{abk}$ into 
$\frac{\sx_{ab} \sx_{ak} }{ \sx_{a'c} \sx_{bk} }$.
Then this expression is pulled through
all $\SF$'s between $\SF_{abk}$ and
$\SF_{a'ck}$. They act only on the first indices of
$\sx_{ak}$ and $\sx_{bk}$ as the following chain
of permutations:
$$
P_{a+1,b}\ldots P_{b-1,b} \ P_{1,b+1} \ldots P_{b,b+1} 
 \ldots  P_{1,c} \ldots P_{a'-1,c} \,.
$$
The part  $P_{a+1,b}\ldots P_{b-1,b}$ transforms 
$\frac{ \sx_{ak} }{ \sx_{bk} }$ into
$\frac{ \sx_{ak} }{ \sx_{a+1,k} }$. 
Then each of the shift operators
$U_d$,\  $b \leq d \leq c-2$ increases the first 
indices of $\sx_{ak}$ and $\sx_{a+1,k}$
by one, which yields  $\frac{ \sx_{a'-1,k} }{ \sx_{a'k} }$.
The last part, $P_{1,c} \ldots P_{a'-1,c}$, transforms
$\frac{ \sx_{a'-1,k} }{ \sx_{a'k} }$ into 
$\frac{ \sx_{ck} }{ \sx_{a'k} }$. 
Finally, $\frac{\sx_{ab} \sx_{ck} }{ \sx_{a'c} \sx_{a'k} }$ 
is transformed by $\SF_{a'ck}$ into 
$\frac{\sx_{ab} }{ \sx_{a'c} }$ which then
is not affected by the remaining~$\SF$'s. 

Thus, we have verified {\bf i)} and {\bf ii)} 
and therefore proved relation~(\ref{RFF}). 

In order to establish Eq.~(\ref{RFFr}) we will
prove first that 
\begin{equation}\label{RFFFr} 
   \phi_0 \bigl(\bar{W}'(3,N)\bigr) = 
   \varphi'' \bigl(\bar{W}'(3,N)\bigr) \cdot \SF_N  \,,
\end{equation}
where $\varphi''$ is a homomorphism from 
${\mathcal B}(3,N)$ to $\text{End}\,\bigl(S_N\bigr)$,
\begin{equation}\label{phi3}   
 \varphi'' ( R_{abc} ) =
 \dil{ \frac{q \sx_{ab} }{ \sx_{a+c-b,c} }
 \prod_{k=0}^{c-b-1} \frac{ \sy_{a+k,b+k}}{\sy_{a+k,b+k+1}} }.  
\end{equation}
Indeed, we have 
\begin{equation}\label{Rhatabc}    
 \SR_{abc} = \phi_0(R_{abc}) = 
 \hat{\SR}_{abc} \cdot \SF_{abc} \,, \qquad 
 \hat{\SR}_{abc} =
 \dil{ q \sx_{ab}   \sx_{bc}^{-1} 
    \sy_{ab} \sy_{ac}^{-1} }   .
\end{equation}
Note that the $\sx$'s arguments of 
$\hat{\SR}_{abc}$ and $\check{\SR}_{abc}$ are the same.
Therefore, moving all $\SF$'s in 
$\phi_0 \bigl(\bar{W}'(3,N)\bigr) $ 
to the right results in the same transformation
of these arguments as moving all $\SF$'s in  
$\phi_0 \bigl(W'(3,N)\bigr)$ to the left.

What the $\sy$'s arguments of $\hat{\SR}_{abc}$
are concerned, they are transformed only by
those $\SF$'s that enter the 
$\phi_0\bigl(\lfloor \bar{W}'(2,c-1)\rfloor_c\bigr)$ 
part. Specifically, $\sy_{ab}$ does not change
whereas $\sy_{ac}$ is transformed by
consecutive action of pairs
$\SF_{a+k,b+1+k,c}\SF_{a+1+k,b+1+k,c}$,
$0 \leq k \leq c-b-1$, which yields the
product of $\sy$'s in~(\ref{phi3}).
 
It remains to pull $\SF_N$ in (\ref{RFFFr})
to the left. Note that $\SF_N^2 =1$.
Using (\ref{Fxx}) and (\ref{Fyy}),
one can verify that
\begin{align}
\label{Frx}
{}&  \SF_N \, \sx_{ab} \, \SF_N 
 = \sx_{b-a,b} \, \prod_{k=1}^{N-b} 
  \frac{\sx_{b-a,b+k}}{\sx_{b-a+1,b+k}} \,, \qquad
    \SF_N \, \sy_{ab} \, \SF_N 
= \prod_{k=1}^{b-a} 
  \frac{\sy_{k,a+k}}{\sy_{k,a+k-1}} \,.
\end{align}
Whence
\begin{align}
\label{Frxx}
{}&   \frac{\sx_{ab}}{\sx_{a+1,b+1}}  \, \SF_N 
 = \SF_N \, \frac{\sx_{b-a,b}}{\sx_{b-a+1,b+1}}  \,,\qquad
     \frac{\sy_{ab}}{\sy_{a,b+1}}  \, \SF_N 
= \SF_N \, \frac{\sy_{b-a+1,b}}{\sy_{b-a+1,b+1}}  .
\end{align} 
Applying these relations to (\ref{phi3}),
we obtain~(\ref{phihat}).
\qed\vspace{2mm}

\noindent {\bf Proof of Lemma~\ref{RRR'}.}\ 
Using (\ref{Zkdef}), (\ref{RabZ}), and (\ref{Zxy}),
we find
\begin{align} 
\label{tauphi1}
 \tau \bigl( \ddil{ abc}  \bigr) =
 \varphi \bigl(  R_{abc} \bigr) \,, \qquad
  \tau \bigl( \mu_1 \bigl( \ddil{abc} \bigr) \bigr) 
  = \varphi' \bigl(  R_{abc} \bigr) .
\end{align} 
Therefore, invoking Lemma~\ref{TW23}, we conclude that
\begin{align}
\label{tauphi2}
{}& \tau \bigl( \ST_N \bigr) =
 \tau \Bigl( \bigl( W'(3,N+1)\bigr)^\star \Bigr)
 = \varphi \bigl( W'(3,N+1)\bigr) \,, \\
{}& \tau \Bigl( \mu_1 \bigl( \ST_N \bigr) \Bigr) = 
 \tau \Bigl( \mu_1 \bigl( \bigl( W'(3,N+1)\bigr)^\star \bigr) \Bigr)
 = \varphi' \bigl(  \bar{W}'(3,N+1) \bigr) .
\end{align} 
These relations along with Lemma~\ref{RR'} 
yield the statement of Lemma~\ref{RRR'}.
\qed\vspace*{2mm}

{\bf Acknowledgements.} 
This work was supported in part
by the Swiss National Science Foundation grants
200020-129609 (A.V.) and 200020-126817 (A.B.) 
and by the RFBR grant  11-01-00570.
The authors thank Anton Alekseev, 
Ludwig Faddeev, Sergey Fomin, and Rinat Kashaev 
for helpful conversations.

{\sc \small
\noindent
Section of Mathematics, University of Geneva, 
% 2--4 rue du Li\`evre, 
C.P. 64, 1211 Gen\`eve 4, Switzerland \\
Steklov Mathematical Institute,
Fontanka 27, 191023, St. Petersburg, Russia}
  
\end{document}